\documentclass[a4paper, onecolumn,  superscriptaddress, 10pt, accepted=2025-12-10, issue=6, volume=7,shorttitle=papers/compositionality-7-6]{compositionalityarticle}
\pdfoutput=1
\usepackage[utf8]{inputenc}
\usepackage[english]{babel}
\usepackage[T1]{fontenc}
\usepackage{float}
\usepackage{amsmath,mathtools,amssymb,abraces,enumitem,amsthm,thmtools}
\usepackage{tikz,tikz-cd}
\usepackage{stmaryrd,accents,microtype}
\definecolor{dark-red}{rgb}{0.5,0.15,0.15}
\usepackage[numbers,sort&compress]{natbib}
\usepackage{hyperref}

\newcommand{\C}{\mathcal{C}}

\newcommand{\de}{\partial}
\newcommand{\p}{\times}
\renewcommand{\vec}{\overrightarrow}
\renewcommand{\P}{\mathbb{P}}

\newtheorem*{thmN}{Theorem}

\newtheorem{thm}{Theorem}[section]
\newtheorem{prop}[thm]{Proposition}

\newtheorem{cor}[thm]{Corollary}

\newtheorem{defnot}[thm]{Definition and notation}
\newcommand{\bdn}{\begin{defnot}}
\newcommand{\edn}{\end{defnot}}
\newcommand{\bp}{\begin{prop}}
\newcommand{\ep}{\end{prop}}
\newcommand{\bth}{\begin{thm}}
\renewcommand{\eth}{\end{thm}}
\newcommand{\bpf}{\begin{proof}}
\newcommand{\epf}{\end{proof}}
\newcommand{\bc}{\begin{cor}}
\newcommand{\ec}{\end{cor}}
\newtheorem{defn}[thm]{Definition}

\newcommand{\bd}{\begin{defn}}
\newcommand{\ed}{\end{defn}}
\newtheorem{nota}[thm]{Notation}

\renewcommand{\top}{{\mathbf{Top}}}
\newcommand{\iso}{\cong}

\renewcommand{\leq}{\leqslant}
\renewcommand{\geq}{\geqslant}

\newcommand{\brm}[1]{{\rm{\mathbf{#1}}}}
\newcommand{\dtop}{{\brm{Flow}}}

\newcommand{\set}{{\brm{Set}}}
\newcommand{\poset}{{\brm{PoSet}^+}}

\DeclareMathOperator{\id}{Id}

\DeclareMathOperator{\Ch}{Ch}

\newcommand{\liminj}{\varinjlim}

\newcommand{\rest}{\!\upharpoonright\!}

\makeatletter
\def\varholim@#1#2{%
	\vtop{\m@th\ialign{##\cr
			\hfil$#1\operator@font holim$\hfil\cr
			\noalign{\nointerlineskip\kern1.5\ex@}#2\cr
			\noalign{\nointerlineskip\kern-\ex@}\cr}}%
}
\def\holimproj{%
	\mathop{\mathpalette\varholim@{\leftarrowfill@\textstyle}}\nmlimits@
}
\def\holiminj{%
	\mathop{\mathpalette\varholim@{\rightarrowfill@\textstyle}}\nmlimits@
}
\makeatother

\DeclareMathOperator{\cosk}{cosk}
\DeclareMathOperator{\COSK}{\vec{\cosk}}

\newcommand{\ddownarrow}{{\downarrow}}

\DeclareMathOperator{\seq}{Seq}
\DeclareMathOperator{\vtx}{Vert}

\newcommand{\vd}{\vec{d}\!_1}

\DeclareMathOperator{\TT}{T}
\newcommand{\LMet}{\mathbf{LvMet}}

\allowdisplaybreaks

\begin{document}
	
\title{Towards a theory of natural directed paths}
\date{}
\author{Philippe Gaucher}
\email{gaucher@irif.fr}
\orcid{0000-0003-0287-6252}
\affiliation{Universit\'e Paris Cit\'e, CNRS, IRIF, F-75013, Paris, France}

\newcommand{\authorsforheader}{P.~Gaucher}
\newcommand{\paperdoi}{https://doi.org/10.46298/compositionality-7-6}
\newcommand{\receiveddate}{2024-03-21}
\newcommand{\accepteddate}{2025-12-10}

\begin{abstract}
	We introduce the abstract setting of presheaf category on a thick category of cubes. Precubical sets, symmetric transverse sets, symmetric precubical sets and the new category of (non-symmetric) transverse sets are examples of this structure. All these presheaf categories share the same metric and homotopical properties from a directed homotopy point of view. This enables us to extend Raussen's notion of natural $d$-path for each of them. Finally, we adapt Ziemia\'{n}ski's notion of cube chain to this abstract setting and we prove that it has the expected behavior on precubical sets. As an application, we verify that the formalization of the parallel composition with synchronization of process algebra using the coskeleton functor of the category of symmetric transverse sets has a category of cube chains with the correct homotopy type.
\end{abstract}
	
\maketitle
\tableofcontents

\section{Introduction}

\subsection*{Presentation} Precubical sets are \textit{de facto} the standard geometric model for directed homotopy for concurrency \cite{DAT_book}. In fact, most of them are even non-positively curved in the sense of \cite[Definition~1.28 and Proposition~1.29]{zbMATH07226006}, or in the worst case proper in the sense of \cite[page~499]{MR3722069}.  The motivation for introducing symmetric transverse sets in \cite{symcub} is to formalize the parallel product with synchronization for process algebra using the associated coskeleton functor \cite[Theorem~4.1.8]{symcub}. Indeed, it is impossible to use the coskeleton functor associated with the category of precubical sets because of its pathological behavior (see \cite[Proposition~3.15]{ccsprecub} and \cite[Definition~3.1.3]{symcub}). However, precubical sets still remain sufficient to model this parallel product by tweaking the coskeleton functor of this category (see \cite[Section~3.3]{ccsprecub}).

Symmetric transverse sets share with precubical sets similar metric and homotopical properties by \cite{DirectedDegeneracy}. Indeed, their geometric realization carries a Lawvere metric structure which enables us to extend Raussen's notion of (tame) natural $d$-path originally defined for precubical sets \cite[Definition~2.14]{MR2521708} \cite[Definition~5.3]{MR3722069} \cite[Section~2.9]{MR4070250}. Moreover, the full subcategory of representable objects of this presheaf category is c-Reedy in the sense of \cite[Definition~8.25]{c-Reedy}, like the box category (see Definition~\ref{box_cat}) for the precubical sets (the box category is even direct Reedy in the sense of \cite[Definition~15.1.2]{ref_model2}). This makes possible to compare in \cite[Theorem~7.4]{DirectedDegeneracy} the natural realization of a symmetric transverse set with other realization functors and to generalize homotopical results proved in \cite{realization} and \cite{NaturalRealization} for precubical sets. 

The technical contribution of this note is threefold. Firstly, we explain why precubical sets and symmetric transverse sets belong to a larger family of presheaf categories on a \textit{thick category of cubes} (see Definition~\ref{thick}). This family of presheaf categories contains also the symmetric precubical sets of \cite{MR1988396} and a new category of non-symmetric transverse sets. Symmetric transverse sets are presheaves on a thick category of cubes which turns out to be, for tautological reasons, the greatest one for the inclusion. 

\begin{thmN} (Proposition~\ref{Orem}, Theorem~\ref{thick-cube} and Theorem~\ref{greatest})
There exists a greatest thick category of cubes for the inclusion not containing the symmetry maps.
\end{thmN}

Secondly, we prove that all results of \cite{DirectedDegeneracy} are valid for all presheaf categories on a thick category of cubes. 

\begin{thmN} (Section~\ref{sec-3} and more specifically Theorem~\ref{same-rea})
All metric and homotopical results of \cite{DirectedDegeneracy} are valid for the category ${\mathcal{A}}^{op}\set$ of $\mathcal{A}$-sets when $\mathcal{A}$ is a thick category of cubes. 
\end{thmN}

Thirdly, we obtain a statement which coincides with (a part of) \cite[Theorem~7.5]{MR4070250} when $\mathcal{A}$ is the box category used to define the precubical sets:

\begin{thmN} (Corollary~\ref{Zgen})
	Let $K$ be a precubical set.  Let $\mathcal{A}$ be a thick category of cubes. The space of tame natural $d$-paths of the free $\mathcal{A}$-set $\mathcal{L}_\mathcal{A}(K)$ generated by $K$ is homotopy equivalent to the classifying space of the small category of Ziemia\'{n}ski's cube chains of the free $\mathcal{A}$-set $\mathcal{L}_{\mathcal{A}}(K)$ generated by $K$.
\end{thmN}

This leads to the following application:

\begin{thmN} (Section~\ref{app})
	The formalization of the parallel composition with synchronisation of \cite{symcub} using the coskeleton functor of the category of symmetric transverse sets has a category of cube chains which gives the correct space of tame natural $d$-paths up to homotopy.
\end{thmN}

The formal setting of presheaf category on a thick category of cubes is a first step towards an axiomatization of the notion of tame natural $d$-path. The next step would be to find a way of taking into account the globular version of this notion as it is introduced in \cite{Moore3}. This could lead to a general framework unifying all geometric approaches of directed homotopy for concurrency~\footnote{The expression ``directed homotopy'' has several quite distinct meanings. It is the reason why I add ``for concurrency'' on purpose.}.

\subsection*{Acknowledgments} I am grateful to the anonymous referee for pointing out the flaw in the proof of Theorem~\ref{cube-chain-morphism} (which requires to add Proposition~\ref{newmetric}) and for the suggestions to improve the presentation of the paper.

\section{Thick category of cubes}

\bd \cite{LawvereMetric} \label{Lawvere_def}
The small monoidal category $([0,\infty],\geq,+,0)$ has for objects the interval $[0,\infty]$ and there is a unique arrow $x\to y$ whenever $x\geq y$. It is equipped with the monoidal structure induced by addition. A small category enriched over $([0,\infty],\geq,+,0)$ is called a \textit{Lawvere metric space}. The category of Lawvere metric spaces is denoted by $\LMet$.
\ed

Let us now expand Definition~\ref{Lawvere_def}. A \textit{Lawvere metric space} $(X,d)$ is a set $X$ equipped with a map $d:X\p X\to [0,\infty]$ called a \textit{Lawvere metric} such that:
\begin{itemize}
	\item $\forall x\in X,d(x,x)=0$
	\item $\forall (x,y,z)\in X\p X\p X, d(x,y)\leq d(x,z)+d(z,y)$.
\end{itemize}
A map $f:(X,d)\to (Y,d)$ of Lawvere metric spaces is a function $f:X\to Y$ which is \textit{non-expansive}, i.e. $\forall (x,y)\in X\p X, d(f(x),f(y))\leq d(x,y)$.

\begin{nota}
	The category of partially ordered sets or posets together with the \textit{strictly increasing maps} is denoted by $\poset$.  
\end{nota}

Let $[0] = \{()\}$ and $[n] = \{0<1\}^n$ for $n \geq 1$ equipped with the product order. Let $0_n=(0,\dots,0)$ ($n$ times) and $1_n=(1,\dots,1)$ ($n$ times) with $n\geq 0$. By convention, one has $\{0<1\}^0=[0]=\{()\}$. In the sequel, for all $n\geq 1$, both the sets $[n]$ and $[0,1]^n$ are equipped with the product order. By convention, $[0,1]^0$ is a singleton.

Let $\delta_i^\alpha : [n-1] \rightarrow [n]$ be the \textit{coface map} defined for $1\leq i\leq n$ and $\alpha \in \{0,1\}$ by \[\delta_i^\alpha(x_1, \dots, x_{n-1}) = (x_1,\dots, x_{i-1}, \alpha, x_i, \dots, x_{n-1}).\] 

\bd \label{box_cat} The \textit{box category} $\square$ is the subcategory of $\poset$ generated by the coface maps $\delta_i^\alpha$.  
\ed

Let $x=(x_1,\dots,x_n)$ and $x'=(x'_1,\dots,x'_n)$ be two elements of $[0,1]^n$ with $n\geq 1$. Let $\vd:[0,1]^n\p [0,1]^n \to [0,\infty]$ be the function defined by
\[
\vd(x,x') = 
\begin{cases}
	\displaystyle\sum\limits_{i=1}^{n} |x_i-x'_i| & \hbox{ if } x\leq x'\\
	\infty & \hbox {otherwise.}
\end{cases}
\]
Let $n\geq 0$. The function $\vd:[0,1]^n\p [0,1]^n\to [0,\infty]$ is a Lawvere metric by \cite[Proposition~1.5]{DirectedDegeneracy}. It restricts to a Lawvere metric on $\{0<1\}^n$.

\bd \label{cotransverse} \cite[Definition~2.1.5]{symcub} 
A map $f:[m] \rightarrow [n]$ of $\poset$ is \textit{cotransverse} if \[\hbox{For all } x,y\in [m], \vd(x,y)=1 \hbox{ implies }\vd(f(x),f(y)) = 1.\] Denote by $\widehat{\square}_S$ the subcategory of $\poset$ consisting of the cotransverse maps.
\ed 

A \textit{cotransverse degeneracy map} is a cotransverse map $[n]\to [n]$ for $n\geq 2$ which is not one-to-one. Proposition~\ref{decomposition_distance} and Proposition~\ref{square} are important for the sequel.

\bp \cite[Proposition~3.1.14]{symcub} \label{decomposition_distance} Let $0\leq m\leq n$. Every cotransverse map $f:[m] \rightarrow [n]$ factors uniquely as a composite $[m] \stackrel{\psi}\longrightarrow [m] \stackrel{\phi}\longrightarrow [n]$ with $\phi\in \square$ and $\psi$ cotransverse.  \ep

\bp \label{square}
Let $f:[m]\to [n]$ be a cotransverse map. Let $\delta:[n]\to [p]$ be a map of $\square$. Suppose that $\delta{f}\in \square$. Then $f\in\square$.
\ep

\bpf
There exists a function $s:[p]\to [n]$ obtained by removing some coordinates such that $s\delta=\id_{[n]}$. We deduce that $f = s(\delta{f})$. From $\delta{f}\in \square$, we then deduce that $f\in \square$.
\epf

Let $\sigma_i:[n] \rightarrow [n]$ be the function defined for $1\leq i\leq n-1$ and $n\geq 2$ by \[\sigma_i(x_1, \dots, x_{n}) = (x_1, \dots, x_{i-1},x_{i+1},x_{i}, x_{i+2},\dots,x_{n}).\] These maps are called the \textit{symmetry maps} \cite{MR1988396}. The symmetry maps are clearly cotransverse. This is the reason of the $S$ in the notation $\widehat{\square}_S$.

\begin{nota}
	The subcategory of $\poset$ generated by the coface maps and the symmetry maps is denoted by $\square_S$. 
\end{nota}

\bd \cite[Definition~2.1.7 and Definition~2.1.12]{symcub}
A \textit{category of cubes} is a small category $\mathcal{A}$ satisfying the inclusions \[\square\subset \mathcal{A} \subset \widehat{\square}_S.\] A presheaf on $\mathcal{A}$ is called an \textit{$\mathcal{A}$-set}. The category of $\mathcal{A}$-sets is denoted by $\mathcal{A}^{op}\set$. 
\ed

The category $\square$, $\square_S$ and $\widehat{\square}_S$ are examples of categories of cubes. The $\square$-sets are the \textit{precubical sets} \cite{Brown_cube}. The $\square_S$-sets are the \textit{symmetric precubical sets} (see \cite{MR1988396}). The $\widehat{\square}_S$-sets are the \textit{symmetric transverse sets} introduced in \cite{DirectedDegeneracy}. 

The inclusion of small categories $j_\mathcal{A}:\square\subset \mathcal{A}$ induces by precomposition a forgetful functor \[\omega_\mathcal{A}:\mathcal{A}^{op}\set \longrightarrow \square^{op}\set\] which has a left adjoint \[\mathcal{L}_\mathcal{A}:\square^{op}\set \longrightarrow \mathcal{A}^{op}\set\] given by the left Kan extension along $j_\mathcal{A}$. For a precubical set $K$, the $\mathcal{A}$-set $\mathcal{L}_\mathcal{A}(K)$ is called the \textit{free $\mathcal{A}$-set} generated by $K$. 

\bp \label{LA-faithful}
The functor $\mathcal{L}_\mathcal{A}:\square^{op}\set \to \mathcal{A}^{op}\set$ is faithful.
\ep

\bpf
By \cite[Proposition~2.1.15]{symcub}, the identity map $\id_{\mathcal{L}_\mathcal{A}(K)}$ induces for all precubical sets $K$ a natural inclusion of precubical sets $K\subset \omega_\mathcal{A}\mathcal{L}_\mathcal{A}(K)$. Let $f,g:K\to L$ be two maps of precubical sets such that $\mathcal{L}_\mathcal{A}(f)=\mathcal{L}_\mathcal{A}(g)$. Then $\omega_\mathcal{A}\mathcal{L}_\mathcal{A}(f)=\omega_\mathcal{A}\mathcal{L}_\mathcal{A}(g)$. Thus 
\[f = \omega_\mathcal{A}\mathcal{L}_\mathcal{A}(f)\rest_K = \omega_\mathcal{A}\mathcal{L}_\mathcal{A}(g)\rest_K = g.\]
\epf

Let $K$ be an $\mathcal{A}$-set. The set $K([n])$ is denoted by $K_n$. The vertex of $x(0_n) \in K_0$ is called the \textit{initial state} of the $n$-cube $c$ and the vertex $x(1_n)\in K_0$ is called the \textit{final state} of the $n$-cube $c$. For any map $k:[m] \rightarrow [n]$ of $\mathcal{A}$ and any $\mathcal{A}$-set $K$, denote by $k^*:K_n \rightarrow K_m$ the function induced by $k$. Let $p\geq 0$. The \textit{$p$-cube} $\mathcal{A}[p]$ is by definition the presheaf $\mathcal{A}(-,[p])$. For any $\mathcal{A}$-set $K$, an element $x\in K_n$ corresponds by the Yoneda lemma to a map of $\mathcal{A}$-sets $x:\mathcal{A}[n]\to K$. For any $\mathcal{A}$-set $K$, the data
\[
(K_{\leq n})_p = \begin{cases}
	K_p & \hbox{ if }p\leq n\\
	\varnothing & \hbox{ if }p> n.
\end{cases}
\]
assemble into an $\mathcal{A}$-set denoted by $K_{\leq n}$ because $\mathcal{A}([m],[n]) = \varnothing$ when $m>n$. Let \[\de\mathcal{A}[n] = \mathcal{A}[n]_{\leq n-1}\] for all $n\geq 0$. Let $A={a_1<\dots<a_k} \subset \{1,\dots,n\}$ and $\epsilon\in \{0,1\}$. The \textit{iterated face map} is defined by $\de^\epsilon_A=\de^\epsilon_{a_1} \de^\epsilon_{a_2} \dots \de^\epsilon_{a_k}$ with $\de^\epsilon_{p} = (\delta^\epsilon_{p})^*$.

\bd \label{thick}
A category of cubes $\mathcal{A}$ is \textit{thick} if the factorization of Proposition~\ref{decomposition_distance} is a factorization in $\mathcal{A}$, i.e $f\in \mathcal{A}$ implies $\psi\in \mathcal{A}$.
\ed

The category of cubes $\widehat{\square}_S$ is thick for tautological reasons. The terminology must be understood as follows. A thick category of cubes $\mathcal{A}$ is morally a thick subcategory of the category of cubes $\widehat{\square}_S$: it is an analogy with the notion of thick subcategory of a triangulated category. Theorem~\ref{thick-cube} provides other examples of thick categories of cubes.

\bp \label{explanation-thick}
Let $\mathcal{A}$ be a thick category of cubes. For $n\geq 0$, let \[j^n_\mathcal{A}:(\square\ddownarrow [n]) \stackrel{\subset}\longrightarrow (\mathcal{A}\ddownarrow[n])\] be the functor between comma categories induced by the inclusion $j_\mathcal{A}:\square\subset \mathcal{A}$. Then for all $n\geq 0$ and for all objects $k$ of $(\mathcal{A}\ddownarrow[n])$, the comma category $(k\ddownarrow j^n_\mathcal{A})$ has an initial object.
\ep

\bpf
Let $k:[p]\to [n]$. Using Proposition~\ref{decomposition_distance} and since $\mathcal{A}$ is thick by hypothesis, we obtain the commutative diagram of $\mathcal{A}$
\[
\begin{tikzcd}[row sep=3em, column sep=3em]
{[p]} \arrow[d,"k"'] \arrow[r,"f"] & {[p]} \arrow[d,"\delta\in \square"] \\
{[n}] \arrow[r,equal] & {[n]}
\end{tikzcd}
\]
which is an element of the comma category $(k\ddownarrow j^n_\mathcal{A})$. Consider another element of $(k\ddownarrow j^n_\mathcal{A})$ depicted by the following commutative diagram of solid arrows of $\mathcal{A}$:
\[
\begin{tikzcd}[row sep=3em, column sep=3em]
	{[p]} \arrow[d,equal] \arrow[r,dashed,"g"]  & {[p]} \arrow[d,dashed,"\delta''\in \square"]\\
{[p]} \arrow[d,"k"'] \arrow[r,"f'"]  & {[q]} \arrow[d,"\delta'\in \square"] \\
{[n]} \arrow[r,equal] & {[n]}
\end{tikzcd}
\]
Since $\mathcal{A}$ is thick, write $f'=\delta''g$ with $\delta''\in \square$. We obtain $\delta{f} = k=\delta' f'=\delta'\delta''g$. By uniqueness of the factorization of Proposition~\ref{decomposition_distance}, we obtain $\delta=\delta'\delta''$ and $f=g$. We have obtained the map of $(k\ddownarrow j^n_\mathcal{A})$ 
\[
\begin{tikzcd}[row sep=1em, column sep=1em]
	{[p]}
	\arrow[rd,equal] \arrow[dd,"k"'] \arrow[rr,"f"]  && {[p]}
	\arrow[rd,"\delta''"] \arrow[dd,"\delta",pos=0.8]  &\\
	& {[p]} \arrow[rr,"f'",pos=0.2,crossing over]   && {[q]} \arrow[dd,"\delta'"] \\
	{[n]} \arrow[rd,equal] \arrow[rr,equal]  &&
	{[n]} \arrow[rd,equal]& \\
	& {[n]} \arrow[rr,equal] \arrow[from=uu,"k"',pos=0.15,crossing over] &&
	{[n]}
\end{tikzcd}
\]
Moreover, the map $\delta'':[p]\to [q]$ is unique because it is given by the factorization of $f'$ using Proposition~\ref{decomposition_distance}. Hence the proof is complete.
\epf

\bp \label{free_square}
Let $\mathcal{A}$ be a category of cubes. For all $n\geq 0$, one has the isomorphism of $\mathcal{A}$-sets \[\mathcal{L}_\mathcal{A}(\square[n]) \iso \mathcal{A}[n].\] If moreover $\mathcal{A}$ is thick, then there is the isomorphism of $\mathcal{A}$-sets \[\mathcal{L}_\mathcal{A}(\de\square[n]) \iso \de\mathcal{A}[n]\] for all $n\geq 0$. 
\ep

\bpf
The first statement is \cite[Proposition 2.1.14]{symcub}. Let $n\geq 0$. Since $\mathcal{L}_\mathcal{A}$ is colimit-preserving, there is a natural map of $\mathcal{A}$-sets ($\square_{<n}$ and $\mathcal{A}_{<n}$ are the full subcategory of $\square$ and $\mathcal{A}$ respectively containing only $[0],\dots,[n-1]$): 
\[
\mathcal{L}_\mathcal{A}(\de\square[n]) \iso \liminj_{(\square_{<n}\ddownarrow [n])} \mathcal{L}_\mathcal{A}(\square[p]) \iso \liminj_{(\square_{<n}\ddownarrow [n])} \mathcal{A}[p] \to \liminj_{(\mathcal{A}_{<n}\ddownarrow [n])} \mathcal{A}[p] \iso \de\mathcal{A}[n].
\]
The above arrow is an isomorphism by Proposition~\ref{explanation-thick} and \cite[Theorem~1 p.~213]{MR1712872}. 
\epf

\bp \label{Orem}
The set of maps \[\widehat{\square} = \{\phi:[m]\to[n]\in \widehat{\square}_S\mid \forall \delta:[p]\to [m]\in \square, \phi\delta \hbox{ one-to-one }\Rightarrow \phi\delta\in \square\}\] is closed under composition and contains all identities of $\widehat{\square}_S$. There are the inclusions $\square\subset \widehat{\square} \subset \widehat{\square}_S$. In other terms, the set of maps $\widehat{\square}$ yields a well-defined category of cubes. The only one-to-one functions of $\widehat{\square}$ are the maps of $\square$. In particular, the only bijective map of $\widehat{\square}([n],[n])$ is $\id_{[n]}$ for all $n\geq 0$: $\widehat{\square}$ does not contain any symmetry map.
\ep

\bpf
Let $\phi_1,\phi_2\in \widehat{\square}$ such that $\phi_1\phi_2$ exists. Let $\delta\in \square$ such that $\phi_1\phi_2\delta$ exists and is one-to-one. Then $\phi_2\delta$ is a one-to-one function. Thus $\phi_2\delta\in \square$, $\phi_2$ belonging to $\widehat{\square}$. We deduce that $\phi_1\phi_2\delta =\phi_1(\phi_2\delta) \in \square$ since $\phi_1\in \widehat{\square}$. This means that $\phi_1\phi_2 \in \widehat{\square}$. For all $\phi=\id$, one has $\phi\delta=\delta \in \square$. Hence $\widehat{\square}$ contains all identity maps. Finally suppose that $f:[m]\to [n]\in \widehat{\square}$ is one-to-one. Then $f\id_{[m]}$ is one-to-one, which implies that $f\in \square$.
\epf

\bd The $\widehat{\square}$-sets are called \textit{transverse sets}. 
\ed

The following maps, introduced in \cite[Definition~3.1.11]{symcub}, are examples of cotransverse degeneracy maps. Let $\gamma_i:[n] \rightarrow [n]$ be the function defined for $1\leq i\leq n-1$ and $n\geq 2$ by \[\gamma_i(x_1, \dots, x_{n}) = (x_1, \dots, x_{i-1},\max(x_{i},x_{i+1}),\min(x_{i},x_{i+1}), x_{i+2},\dots,x_{n}).\] 

\begin{nota} \cite[Theorem~3.1.16]{symcub} \label{smallercube}
	The category of cubes generated by the $\delta_i^\alpha$, $\sigma_i$ and $\gamma_i$ is denoted by $\overline{\square}$. 
\end{nota}

Unlike $\widehat{\square}$, the category of cubes $\overline{\square}$ has a conjectural presentation by generators and relations \cite[Proposition~3.1.20 and Conjecture~3.1.21]{symcub}.

\bth \label{thick-cube}
The categories of cubes $\square$, $\square_S$, $\widehat{\square}$ and $\overline{\square}$ are thick.
\eth

\bpf
Every map $f:[m]\to [n]$ of $\square$ factors uniquely in $\widehat{\square}_S$ as a composite $[m] \to [m] \to [n]$ such that the right-hand map $[m]\to [n]$ belongs to $\square$. Since $f=f\id_{[m]}$, we deduce by uniqueness that the left-hand map belongs to $\square$. Hence the category $\square$ is thick.

Every map $f:[m]\to [n]$ of $\square_S$ factors uniquely in $\widehat{\square}_S$ as a composite $[m] \to [m] \to [n]$ such that the right-hand map $[m]\to [n]$ belongs to $\square$. Since all maps of $\square_S$ are one-to-one, the left-hand map $[m] \to [m]$ is one-to-one, hence bijective for cardinality reason. We deduce that the left-hand map is a composite of symmetry maps, which means that it belongs to $\square_S$. Thus the category $\square_S$ is thick. 

Consider a map $f:[m]\to [n]$ of $\widehat{\square}$. It  factors in $\widehat{\square}_S$ as a composite $[m]\to [m]\to [n]$ such that the right-hand map $\delta:[m]\to [n]$ belongs to the box category. Denote by $g:[m]\to [m]$ the left-hand map. Let $\delta':[p]\to [m]\in \square$ such that $g\delta'$ exists and is one-to-one. Then $f\delta'=\delta(g\delta')$ is one-to-one, being the composite of two one-to-one functions. Since $f\in \widehat{\square}$, we deduce that $\delta{g}\delta'\in \square$. Using Proposition~\ref{square}, we deduce that $g\delta'\in\square$. We have proved that $g\in \widehat{\square}$. This means that the category of cubes $\widehat{\square}$ is thick. 

Finally, one has $\gamma_i\delta_j^\alpha =\delta_j^\alpha \gamma_i$ for $j\leq i-1$ and $j\geq i+2$, $\gamma_i\delta_i^0=\delta_{i+1}^0$, $\gamma_i\delta_i^1=\delta_i^1$, $\gamma_i\delta_{i+1}^0=\delta_{i+1}^0$ and $\gamma_i\delta_{i+1}^1=\delta_{i}^1$. Thus, $\overline{\square}$ is thick.
\epf

It is possible to obtain even more examples of thick categories of cubes by removing the symmetry maps from $\overline{\square}$, and also by adding one by one other cotransverse degeneracy maps in the set of generators. The only interest of the category of cubes $\overline{\square}$ for this note is to show that it is really easy to construct other thick categories of cubes. Theorem~\ref{greatest} gives a characterization of $\widehat{\square}$.

\bth \label{greatest}
Let $\mathcal{A}$ be a thick category of cubes which does not contain the symmetry maps, i.e. the only bijective map of $\mathcal{A}([n],[n])$ is $\id_{[n]}$ for all $n\geq 0$. Then $\mathcal{A}\subset \widehat{\square}$. In other terms, $\widehat{\square}$ is the greatest thick category of cubes for the inclusion which does not contain any symmetry map.
\eth

\bpf
Let $\phi:[m]\to [n]$ be a map of $\mathcal{A}$. Let $\delta:[p]\to [m]\in \square$ such that $\phi\delta$ is one-to-one. Write $\phi\delta=\delta'\phi'$ with $\delta':[p]\to [n]\in \square$ and $\phi':[p]\to [p]\in \mathcal{A}$, the category of cubes $\mathcal{A}$ being thick by hypothesis.  Since $\phi\delta$ is one-to-one, the function $\phi':[p]\to [p]$ is one-to-one. By hypothesis, this implies that $\phi'$ is the identity map. Thus $\phi\delta=\delta'\in \square$. We have proved that $\phi\in \widehat{\square}$.
\epf

To summarize, $\widehat{\square}$ is the greatest thick category of cubes for the inclusion not containing the symmetry maps and $\widehat{\square}_S$ is the greatest thick category of cubes for the inclusion. 

For the sequel, \textit{$\mathcal{A}$ denotes a fixed thick category of cubes}. Note that with the choice $\mathcal{A}=\square$ (the least thick category of cubes), all following results remain valid. However, their formulation is not necessarily the best one.

\section{Metric and homotopical study}
\label{sec-3}

For the ease of the reader, we recall a few definitions from \cite{c-Reedy} before generalizing the results of \cite{DirectedDegeneracy}. The point is not to overload this section but to make it comprehensible. 

\bd \cite[Definition~6.12]{c-Reedy} 
Let $\C$ be a category equipped with an ordinal degree function on its objects. 
\begin{itemize}
	\item A morphism is \textit{level} if its domain and codomain have the same degree.
	\item The \textit{degree} of a factorization $(h,g)$ of a morphism $f$ is the degree of the intermediate object (i.e. the domain of $h$ which is the codomain of $g$).
	\item A factorization of a morphism $f$ is \textit{fundamental} if its degree is strictly less than the degrees of both the domain and codomain of $f$.
	\item A morphism is \textit{basic} if it does not admit any fundamental factorization. 
\end{itemize}
\ed

\bd \cite[page~37]{c-Reedy}
Let $\C$ be a category equipped with an ordinal degree function on its objects. The \textit{$\delta$-th stratum} of $\C$, denoted by $\C_{=\delta}$, is the subcategory of $\C$ generated by the objects of degree $\delta$ and by the basic morphisms between them.
\ed

\bd
Let $\C$ be a category. Let $f$ be a map of $\C$. The \textit{category of factorizations} of $f$ has for objects the pairs $(h,g)$ such that $hg=f$ and for morphisms $k:(h,g)\to (h',g')$ the morphisms $k$ of $\C$ (which are called \textit{connecting morphisms}) such that there is a commutative diagram 
\[
\begin{tikzcd}[row sep=3em,column sep=5em]
	\bullet \arrow[r,"g'"] \arrow[equal,d]  & \bullet \arrow[r,"h'"]& \bullet \arrow[equal,d] \\
\bullet \arrow[r,"g"] & \arrow[u,"k"]\bullet \arrow[r,"h"] & \bullet 
\end{tikzcd}
\]
\ed

\bd \cite[Definition~8.25]{c-Reedy} 
A \textit{c-Reedy category} $\C$ is a small category equipped with an ordinal degree function $d$ on its objects, and subcategories $\overleftrightarrow{\C}$, $\overrightarrow{\C}$ and $\overleftarrow{\C}$ containing all objects such that 
\begin{enumerate}
	\item $\overleftrightarrow{\C}\subseteq \overrightarrow{\C} \cap\overleftarrow{\C}$.
	\item Every morphism in $\overleftrightarrow{\C}$ is level.
	\item Every morphism in $\overrightarrow{\C}\backslash \overleftrightarrow{\C}$ strictly raises degree, and every morphism in $\overleftarrow{\C}\backslash \overleftrightarrow{\C}$ strictly lowers degree.
	\item Every morphism $f$ factors as $\overrightarrow{f}\overleftarrow{f}$, where $\overrightarrow{f}\in \overrightarrow{\C}$ and $\overleftarrow{f}\in \overleftarrow{\C}$. The subcategory of the category of factorizations of $f$ generated by the pairs $(h,g)$ with $h\in\overrightarrow{\C}$ and $g\in \overleftarrow{\C}$ and such that the connecting morphisms belong to $\overleftrightarrow{\C}$ is connected for all $f$. 
	\item For any object $x$ and any degree $\delta<d(x)$, the functor $\overleftarrow{\C}(x,-):\overleftrightarrow{\C}_{=\delta} \to \set$ is a coproduct of retracts of representables.
\end{enumerate}
\ed

\begin{nota}
	Let
	\begin{align*}
		& \vec{\mathcal{A}}=\mathcal{A},\\
		& \overleftrightarrow{\mathcal{A}}=\overleftarrow{\mathcal{A}}= \coprod_{n\geq 0}\{f:[n]\to [n]\mid f\in \mathcal{A}\}.
	\end{align*}
	We consider the degree function $d([n])=n$ for all $n\geq 0$.
\end{nota}

Every morphism $f:[m]\to [n]$ of $\mathcal{A}$ is basic since every factorization of $f$ as a composite $[m]\to [p]\to [n]$ implies that $m\leq p\leq n$, and therefore that every factorization is not fundamental: $p<\min(m,n) =m$ is impossible indeed. Hence, for all $n\geq 0$, the $n$-th stratum $\mathcal{A}_{=n}$ is the full subcategory of $\mathcal{A}$ having one object $[n]$. Moreover, one has \[\mathcal{A}_{=n}([n],[n])=\mathcal{A}([n],[n])\] for all $n\geq 0$.

\bp \label{example-c-Reedy}
 The small category $\mathcal{A}$ is c-Reedy.
\ep

\bpf
Let $f$ be a map of $\mathcal{A}$. Consider the subcategory $\mathcal{F}_f$ of the category of factorizations of $f$ generated by the pairs $(h,g)$ such that $h\in\overrightarrow{\mathcal{A}}$ and $g\in \overleftarrow{\mathcal{A}}$. Note that the connecting maps necessarily belong to $\overleftrightarrow{\mathcal{A}}$ since $\overleftrightarrow{\mathcal{A}}=\overleftarrow{\mathcal{A}}$. By \cite[Proposition~5.8]{DirectedDegeneracy} and since the factorization of Proposition~\ref{decomposition_distance} restricts to a factorization in $\mathcal{A}$, the category $\mathcal{F}_f$ has a final object. The rest of the proof goes like the proof of \cite[Proposition~5.9]{DirectedDegeneracy}. Let us summarize the argument. One has $\overleftrightarrow{\mathcal{A}}\subset \vec{\mathcal{A}} \cap \overleftarrow{\mathcal{A}}$ (first axiom). Every morphism of $\overleftrightarrow{\mathcal{A}}$ is degree-preserving (second axiom). Every morphism of $\vec{\mathcal{A}}\backslash \overleftrightarrow{\mathcal{A}}$ strictly raises degree and every morphism of $\overleftarrow{\mathcal{A}}\backslash \overleftrightarrow{\mathcal{A}}=\varnothing$ strictly lowers degree (third axiom). The category $\mathcal{F}_f$ is connected since it has a final object (fourth axiom). For every $n\geq 0$, and any degree $m<n$, the functor $\overleftarrow{\mathcal{A}}([n],-):\mathcal{A}_{=m}\to \set$ is an (empty) coproduct of retracts of representables because $\mathcal{A}([n],[m])=\varnothing$ (fifth axiom).
\epf

\begin{nota}
	Let $\C$ be a small category. Let $\mathcal{M}$ be a locally small category. The category of functors from $\C$ to $\mathcal{M}$ together with the natural transformations is denoted by $\mathcal{M}^\C$.
\end{nota}

\begin{nota}
	Let $n\geq 0$. Following the notations of \cite[page~37]{c-Reedy}, let 
	\[
	\de_n\mathcal{A}([p],[q]) = \int^{[m]\in \mathcal{A}_{<n}} \mathcal{A}([m],[q]) \p \mathcal{A}([p],[m])
	\]
	The latching and matching object functors $L_n,M_n:\mathcal{M}^{\mathcal{A}}\to \mathcal{M}^{\mathcal{A}_{=n}}$ are given by
	\begin{align*}
		& (M_nA)_{[n]} = \int_{[m]\in \mathcal{A}} A([m])^{\de_n\mathcal{A}([n],[m])} \\
		& (L_nA)_{[n]} = \int^{[p]\in \mathcal{A}}\de_n\mathcal{A}([p],[n]) . A([p]) 
	\end{align*}
\end{nota}

We obtain: 

\bth \label{c-Reedy-model}
Let $\mathcal{M}$ be a model category. Suppose that the projective model structure on $\mathcal{M}^{\mathcal{A}_{=n}}$ exists for all $n\geq 0$. There exists a unique model structure on $\mathcal{M}^{\mathcal{A}}$ such that 
\begin{itemize}
	\item The weak equivalences are objectwise.
	\item A map $A\to B$ of $\mathcal{M}^{\mathcal{A}}$ is a fibration (trivial fibration resp.) if for all $n\geq 0$, the map $A([n])\to (M_nA)_{[n]}\p_{(M_nB)_{[n]}} B([n])$ is a fibration (trivial fibration resp.) of $\mathcal{M}$.
	\item A map $A\to B$ of $\mathcal{M}^{\mathcal{A}}$ is a cofibration (trivial cofibration resp.)if for all $n\geq 0$, $L_nB\sqcup_{L_nA} A\to B$ is a projective cofibration (trivial cofibration resp.) of the projective model structure of $\mathcal{M}^{\mathcal{A}_{=n}}$.
\end{itemize}
This model structure is called the \textit{c-Reedy model structure} of $\mathcal{M}^{\mathcal{A}}$. 
\eth

\bpf
By Proposition~\ref{example-c-Reedy} and \cite[Theorem~8.26]{c-Reedy}, the small category $\mathcal{A}$ is almost c-Reedy in the sense of \cite[Definition~8.8]{c-Reedy}. The proof is complete thanks to \cite[Theorem~8.9]{c-Reedy}.
\epf

\bp\label{de_calcul_func}
One has 
\[
\de_n\mathcal{A}([p],[q]) = \begin{cases}
	\varnothing & \hbox{ if } p>q \hbox{ or }n\leq p\\
	\mathcal{A}([p],[q]) & \hbox{ if } p\leq q \hbox{ and }p<n
\end{cases}
\]
\ep

\bpf
It is mutatis mutandis the proof of \cite[Proposition~5.13]{DirectedDegeneracy}: it suffices to change the category of cubes in the proof and to use Proposition~\ref{decomposition_distance} which restricts to a factorization in $\mathcal{A}$ by definition of a thick category of cubes.
\epf

\bth \label{proj-cof-suff-cond}
Let $\mathcal{M}$ be a model category. Suppose that the projective model structure on $\mathcal{M}^{\mathcal{A}_{=n}}$ exists for all $n\geq 0$. Then the projective model structure on $\mathcal{M}^{\mathcal{A}}$ exists and coincides with the c-Reedy model structure.
\eth

\bpf 
The proof follows the road map of the proof of \cite[Theorem~5.17]{DirectedDegeneracy} and makes use of Theorem~\ref{c-Reedy-model} and Proposition~\ref{de_calcul_func}.
\epf

\begin{nota}
	The category of \textit{$\Delta$-generated spaces} or of \textit{$\Delta$-Hausdorff $\Delta$-generated spaces} (cf. \cite[Section~2 and Appendix~B]{leftproperflow}) is denoted by $\top$. There are three well-known model category structures on it, the $\{q,h,m\}$-model category structures.
\end{nota}

\bd \label{topologize1} \cite[Definition~3.2]{DirectedDegeneracy}
Let $f=(f_1,\dots,f_n):[n]\to [n]$ be a cotransverse map. Let $\TT(f):[0,1]^n\to [0,1]^n$ be the map defined by 
\[
\TT(f)(x_1,\dots,x_n) = (\TT(f)_1(x_1,\dots,x_n),\dots,\TT(f)_n(x_1,\dots,x_n))
\]
with \[\TT(f)_i(x_1,\dots,x_n) = \max_{(\epsilon_1,\dots,\epsilon_n)\in f_i^{-1}(1)} \min \{x_k\mid \epsilon_k=1\}\] for all $1\leq i\leq n$.
\ed

\begin{nota}
	For $\delta_i^\alpha:[n-1]\to [n] \in \square$, let \[\TT(\delta_i^\alpha)=
	\begin{cases}
		[0,1]^{n-1} \to [0,1]^n\\
		(\epsilon_1, \dots, \epsilon_{n-1})\mapsto (\epsilon_1,\dots, \epsilon_{i-1}, \alpha, \epsilon_i, \dots, \epsilon_{n-1})
	\end{cases}
	\] for all $n\geq 1$ and $\alpha\in\{0,1\}$.
\end{nota}

The three mappings $[n]\mapsto [0,1]^n$ for $n\geq 0$, $f:[n]\to [n]\in \mathcal{A} \mapsto \TT(f)$ and $\delta_i^\alpha:[n-1]\to [n] \mapsto \TT(\delta_i^\alpha)$ for $n\geq 1$ give rise to a functor from $\mathcal{A}\subset \widehat{\square}_S$ to $\top$ denoted by $|\mathcal{A}[*]|_{geom}$ by \cite[Theorem~3.9]{DirectedDegeneracy} and to a functor from $\mathcal{A}\subset \widehat{\square}_S$ to $\LMet$ denoted by $|\mathcal{A}[*]|_{\vd}$ by \cite[Theorem~3.16]{DirectedDegeneracy}. Let $K$ be an $\mathcal{A}$-set. Let 
\[
|K|_{geom} = \int^{[n]\in \mathcal{A}} K_n.|\mathcal{A}[n]|_{geom} \hbox{ and }|K|_{\vd} = \int^{[n]\in \mathcal{A}} K_n.|\mathcal{A}[n]|_{\vd}.
\]
These give rise to two colimit-preserving functors from $\mathcal{A}$-sets to topological spaces and Lawvere metric spaces respectively. The latter functor factoring as a composite $\mathcal{A}^{op}\set \to \widehat{\square}_S^{op}\set \to \LMet$, we define the notion of (tame or not) natural $d$-path in the geometric realization $|K|_{geom}$ like in \cite[Section~4]{DirectedDegeneracy}. In fact, by using the inclusion $\mathcal{A}\subset \widehat{\square}_S$ and the fact that $\mathcal{A}$ is thick, we can mimick all constructions of \cite{DirectedDegeneracy} and recover, thanks to Theorem~\ref{proj-cof-suff-cond}, the results already proved for $\widehat{\square}_S$. More precisely, we obtain what follows.

\bd \cite[Definition~4.11]{model3} \label{def-flow}
A \textit{flow} is a small semicategory enriched over the closed monoidal category $(\top,\p)$. The corresponding category is denoted by $\dtop$. The objects are called \textit{states} and the morphisms \textit{execution paths}. The set of states of a flow $X$ is denoted by $X^0$. The space of execution paths from $\alpha$ to $\beta$ of a flow $X$ is denoted by $\P_{\alpha,\beta}X$. 
\ed

There is an inclusion functor $\poset\subset \dtop$ such that there is an execution path from $\alpha$ to $\beta$ if and only if $\alpha<\beta$. 

Let $r\in \{q,m,h\}$. The category $\dtop$ is equipped with its r-model structure \cite[Theorem~7.4]{QHMmodel}. The weak equivalences of the r-model structure are the maps of flows $f:X\to Y$ inducing a bijection on states and a weak equivalence of the r-model structure of $\top$ for each map $\P_{\alpha,\beta}X\to \P_{f(\alpha),f(\beta)}Y$ for $(\alpha,\beta)$ running over $X^0\p X^0$. The fibrations of the r-model structure are the maps of flows $f:X\to Y$ inducing a fibration of the r-model structure of $\top$ for each map $\P_{\alpha,\beta}X\to \P_{f(\alpha),f(\beta)}Y$ for $(\alpha,\beta)$ running over $X^0\p X^0$. The three model category structures of $\dtop$ are accessible in the sense of \cite[Definition~5.1]{zbMATH06722019}. The q-model structure of $\dtop$ is even combinatorial.

\bd \label{def-rea-flow} Let $r\in \{q,m,h\}$. A functor $F:\mathcal{A}^{op}\set \to \dtop$ is a \textit{r-realization functor} if it satisfies the following properties:
\begin{itemize}
	\item $F$ is colimit-preserving.
	\item For all $n\geq 0$, the map $F(\de\mathcal{A}[n])\to F(\mathcal{A}[n])$ is a r-cofibration of $\dtop$.
	\item There is a map $F(\mathcal{A}[*])\to \{0<1\}^*$ in $\dtop^{\mathcal{A}}$ which is an objectwise weak equivalence of the r-model structure of $\dtop$.
\end{itemize}
\ed

\bd \label{cofibrant-A}
An $\mathcal{A}$-set $K$ is \textit{cellular} if the canonical map $\varnothing\to K$ is a transfinite composition of pushouts of the maps $\de\mathcal{A}[n]\to \mathcal{A}[n]$ for $n\geq 0$. An $\mathcal{A}$-set $K$ is \textit{cofibrant} if it is a retract of a cellular $\mathcal{A}$-set. 
\ed

The composite functor $\mathcal{A}\subset \widehat{\square}_S \to \dtop$ taking $[n]$ to the flow $|\square_S[n]|_{nat}$ defined in \cite[Proposition~7.1]{DirectedDegeneracy} induces a colimit-preserving functor 
\[
|K|_{nat} = \int^{[n]\in \mathcal{A}} K_n.|\square_S[n]|_{nat}
\]
from $\mathcal{A}^{op}\set$ to $\dtop$ which is a m-realization functor by \cite[Theorem~7.4]{DirectedDegeneracy}. The composite functor $\mathcal{A}\subset \widehat{\square}_S \to \dtop$ taking $[n]$ to the flow $(\{0<1\}^n)^{cof}$, where $(-)^{cof}$ is a q-cofibrant replacement of $\dtop$, induces a colimit-preserving functor 
\[
|K|_{q} = \int^{[n]\in \mathcal{A}} K_n.(\{0<1\}^n)^{cof}
\]
from $\mathcal{A}^{op}\set$ to $\dtop$ which is a q-realization functor by \cite[Theorem~6.7]{DirectedDegeneracy}.

\bth \label{same-rea} 
There exists an m-realization functor $F:\mathcal{A}^{op}\set \to \dtop$ and two natural transformations inducing bijections on the sets of states \[|-|_q\Longleftarrow F(-)\Longrightarrow |-|_{nat}\] such that for all cofibrant $\mathcal{A}$-sets $K$ and all $(\alpha,\beta)\in K_0\p K_0$, there is the zigzag of natural homotopy equivalences between m-cofibrant topological spaces 
\[\begin{tikzcd}
		\P_{\alpha,\beta} |K|_q & \arrow[l,"\simeq"']  \P_{\alpha,\beta} F(K) \arrow[r,"\simeq"] &\P_{\alpha,\beta} |K|_{nat}.
\end{tikzcd}
\] 
\eth

\bpf
We follow the proof of \cite[Theorem~7.4]{DirectedDegeneracy}. Details are left to the reader.
\epf

\section{Ziemia\'{n}ski's cube chain}

We define at first the category of cube chains of an $\mathcal{A}$-set. We make the link with Ziemia\'{n}ski's original notion in Corollary~\ref{link}.

\begin{nota}
	Let $n\geq 0$. Let $\seq(n)$ be the set of sequences of positive integers $\underline{n}=(n_1,\dots,n_p)$ with $n_1+\dots + n_p=n$. For $\underline{n}\in \seq(n)$, let \[\vtx(\underline{n}) = \bigg\{\sum_{i=1}^{i=j}n_j \mid 0\leq j\leq p\bigg\}\] the set of \textit{vertices} of $\underline{n}$. The number $|\underline{n}|=n$ is the \textit{length} of $\underline{n}$ and $\ell(\underline{n})=p$ is the \textit{number of elements} of $\underline{n}$. 
\end{nota}

\bd 
Let $\underline{n}\in \seq(n)$ with $n\geq 0$. The \textit{$\underline{n}$-cube} is the $\mathcal{A}$-set 
\[
\mathcal{A}[\underline{n}] = \mathcal{A}[n_1] * \dots * \mathcal{A}[n_p]
\]
where the notation $*$ means that the final state $1_{n_i}$ of the $\mathcal{A}$-set $\mathcal{A}[n_i]$ is identified with the initial state $0_{n_{i+1}}$ of the $\mathcal{A}$-set $\mathcal{A}[n_{i+1}]$ for $1\leq i\leq p-1$. 
\ed 

The category of posets being cocomplete, being locally presentable, the vertices of the $\underline{n}$-cube can be equipped with the poset structure induced by the one on each $\mathcal{A}[n_i]_0=\{0<1\}^{n_i}$.

The set $\vtx(\underline{n})$ is in bijection with the set consisting of the initial states of the $\mathcal{A}[n_i]$ for $1\leq i \leq p$ and the final state of $\mathcal{A}[n_p]$ viewed as subobjects of $\mathcal{A}[\underline{n}]$. For all $p,q\in \vtx(\underline{n})$, $p\leq q$ as integers if and only if $p\leq q$ for the poset structure on the vertices of $\mathcal{A}[\underline{n}]$. 

Using the fact that the $\underline{n}$-cube is a colimit, and $\LMet$ being cocomplete, the set of vertices of the $\underline{n}$-cube can also be equipped with the $\vd$ metric. Proposition~\ref{newmetric} summarizes what is necessary to know about this Lawvere metric for the proof of Theorem~\ref{cube-chain-morphism}. 

\bp \label{newmetric}
Let $\underline{n}\in \seq(n)$ with $n\geq 0$ as above. The following properties hold:
\begin{enumerate}[leftmargin=*]
	\item Let $x$ and $y$ be two vertices of $\mathcal{A}[\underline{n}]$ belonging to the same $\mathcal{A}[n_i]$. Then $\vd(x,y)$ is the distance in $\mathcal{A}[n_i]$ from $x$ to $y$.
	\item Let $p,q\in \vtx(\underline{n})$ which are identified with the corresponding vertices of $\mathcal{A}[\underline{n}]$ as explained above. Then $p\leq q$ implies $\vd(p,q)=q-p$ and $p>q$ implies $\vd(p,q)=\infty$.
	\item If $x\in \mathcal{A}[n_{i}]_0$ and $y\in \mathcal{A}[n_{j}]_0$ with $i<j$, then \[\vd(x,y)= \vd(x,n_1+\dots+n_i) + n_{j-1}-n_i + \vd(n_1+\dots+n_{j-1},y)\] where $n_1+\dots+n_i\in \vtx(\underline{n})$ is identified with the final state of $\mathcal{A}[n_i]$ in $\mathcal{A}[\underline{n}]$ and where $n_1+\dots+n_{j-1}\in \vtx(\underline{n})$ is identified with the initial state of $\mathcal{A}[n_j]$ in $\mathcal{A}[\underline{n}]$.
\end{enumerate}
\ep

\bpf
The category $\LMet$ being the category of $([0,\infty],\geq,+,0)$-enriched small categories, we obtain for two vertices $x$ and $y$ of $\mathcal{A}[\underline{n}]$ the general formula 
\[
\vd(x,y) = \min_{\substack{r\geq 1\\(k_1, \dots,k_r)\in \{1,\dots,p\}^r}}\min_{\substack{(x_0,\dots,x_r)\\x=x_0,y=x_r\\\forall i,\{x_i,x_{i+1}\}\subset \mathcal{A}[n_{k_i}]_0}} \vd(x_0,x_1) + \dots + \vd(x_{r-1},x_r).
\]
Let $x$ and $y$ be two vertices of $\mathcal{A}[\underline{n}]$ belonging to the same $\mathcal{A}[n_i]$. Then $\vd(x,y)$ in $\mathcal{A}[n_i]$ is the minimum by the triangular inequality. Hence $(1)$. Let $p,q\in \vtx(\underline{n})$ with $p\leq q$. For an $s$-cube $\mathcal{A}[s]$, $\vd(0_s,1_s)=s$. We deduce that for any $r\geq 1$, any $(k_1, \dots,k_r)\in \{1,\dots,p\}^r$ and any $(p=x_0,\dots,x_r=q)$ such that $\forall i,\{x_i,x_{i+1}\}\subset \mathcal{A}[n_{k_i}]_0$, the sum $\vd(x_0,x_1) + \dots + \vd(x_{r-1},x_r)$ is equal to $q-p$ or to $\infty$. Hence the first part of $(2)$. Assume now that $q>p$. Then in any sum of the type $\vd(x_0,x_1) + \dots + \vd(x_{r-1},x_r)$ above, one of the terms is equal to $\infty$. Hence $(2)$. We deduce $(3)$ thanks to the following calculation:
\begin{align*}
	&\vd(x,y)\\&=\vd(x,n_1+\dots+n_i) + \vd(n_1+\dots+n_i,n_1+\dots+n_{j-1}) + \vd(n_1+\dots+n_{j-1},y)\\
	&=\vd(x,n_1+\dots+n_i) + n_{j-1}-n_i + \vd(n_1+\dots+n_{j-1},y).
\end{align*}
\epf

\bd Let $K$ be an $\mathcal{A}$-set. Let $\alpha,\beta\in K_0$. Let $n\geq 0$. The \textit{small category \[\Ch_{\alpha,\beta}(K,n)\] of cube chains of $K$} is defined as follows. The objects are the maps of $\mathcal{A}$-sets \[\mathcal{A}[\underline{n}] \longrightarrow K\] with $\underline{n}=(n_1,\dots,n_p)$ and $|\underline{n}|=n$ where the initial state of $\mathcal{A}[n_1]$ is mapped to $\alpha$ and the final state of $\mathcal{A}[n_p]$ is mapped to $\beta$. The morphisms are the commutative diagrams of $\mathcal{A}$-sets of the form
\[
\begin{tikzcd}[row sep=3em, column sep=3em]
	\mathcal{A}[\underline{n}_a] \arrow[d] \arrow[r,"a"] & K \arrow[d,equal] \\
\mathcal{A}[\underline{n}_b]  \arrow[r,"b"]  & K
\end{tikzcd}
\]
such that $|\underline{n}_a| = |\underline{n}_b|$ and $\vtx(\underline{n}_b) \subset \vtx(\underline{n}_a)$. Let 
\[
\Ch(K) = \coprod_{(\alpha,\beta)\in K_0\p K_0} \Ch_{\alpha,\beta}(K).
\]
\ed

To give an explicit description of the morphisms in the category of cube chains, we introduce two families of maps of $\mathcal{A}$-sets $\delta_{i,A,B}$ and $\delta_{\underline{f}}$ in what follows. Let $A\sqcup B=\{1,\dots,m_1+m_2\}$ be a partition with the cardinal of $A$ equal to $m_1>0$ and the cardinal of $B$ equal to $m_2>0$. Let \[\phi_{A,B}:\mathcal{A}[m_1]*\mathcal{A}[m_2] \longrightarrow \mathcal{A}[m_1+m_2]\] be the unique map of $\mathcal{A}$-sets such that 
\begin{align*}
	&\phi_{A,B}(\id_{[m_1]}) = \de^0_B(\id_{[m_1+m_2]}),\\
	&\phi_{A,B}(\id_{[m_2]}) = \de^1_A(\id_{[m_1+m_2]}).
\end{align*}
For $i\in \{1,\dots,\ell(\underline{n})\}$ and a partition $A\sqcup B=\{1,\dots,n_i\}$, let \[{\delta_{i,A,B}=\id_{\mathcal{A}[n_1]}*\dots*\id_{\mathcal{A}[n_{i-1}]}*\phi_{A,B}*\id_{\mathcal{A}[n_{i+1}]}*\dots*\id_{\mathcal{A}[n_{\ell(\underline{n})}]}}.\] For $\underline{f}=(f_1,\dots,f_p)$ with $f_i\in\mathcal{A}([n_i],[n_i])$ for $1\leq i \leq p$, let 
\[
{\delta_{\underline{f}}= f_1*\dots *f_p : \mathcal{A}[\underline{n}] \longrightarrow \mathcal{A}[\underline{n}]}.
\]

\bp \label{decomposition_distance_cor}
Let $\mathcal{A}$ be a thick category of cubes. Let $0\leq m\leq n$. Every map of $\mathcal{A}$-sets $f:\mathcal{A}[m]\to \mathcal{A}[n]$ factors uniquely as a composite $f=\mathcal{L}_\mathcal{A}(h)g$ with $g:\mathcal{A}[m]\to \mathcal{A}[m]$ and $h:\square[m]\to \square[n]$.
\ep

\bpf
It is a rephrasing of Proposition~\ref{decomposition_distance} using the Yoneda lemma.
\epf

\bth \label{cube-chain-morphism}
Let $\underline{n_1}\in \seq(n_1)$ and $\underline{n_2}\in \seq(n_2)$. A map of $\mathcal{A}$-sets from $\mathcal{A}[\underline{n_1}]$ to $\mathcal{A}[\underline{n_2}]$ is a composite of maps of the form $\delta_{i,A,B}$ and $\delta_{\underline{f}}$ if and only if $n_1=n_2$ and $\vtx(\underline{n_2})\subset\vtx(\underline{n_1})$. 
\eth

\bpf
The ``only if'' direction is a consequence of the definitions of $\delta_{i,A,B}$ and $\delta_{\underline{f}}$: note that only a map of the form $\delta_{i,A,B}$ removes vertices. Let us treat the ``if'' part. Let $n=n_1=n_2$. Write $\underline{n_2}=(m_1,\dots,m_r)$ with $\ell(\underline{n_2})=r$. Since $\vtx(\underline{n_2})\subset\vtx(\underline{n_1})$ by hypothesis, there exist $\underline{m_j}\in \seq(m_j)$ for $1\leq j\leq r$ such that $\underline{n_1}=(\underline{m_1},\dots,\underline{m_r})$. The smallest (biggest resp.) element of $\vtx(\underline{n_2})\subset\vtx(\underline{n_1})$ is $0$ ($n$ resp.) and they can be identified with the initial and final states respectively of $\mathcal{A}[\underline{n_1}]$ and $\mathcal{A}[\underline{n_2}]$. Consider a map of $\mathcal{A}$-sets \[k:\mathcal{A}[\underline{n_1}] \longrightarrow \mathcal{A}[\underline{n_2}].\]  Consider two vertices $x,y$ of $\mathcal{A}[\underline{n_1}]$ such that $\vd(x,y)=1$. Then, by Proposition~\ref{newmetric}, the two vertices $x$ and $y$ are the extremities of a unique $1$-cube of $\mathcal{A}[\underline{n_1}]$ going from $x$ to $y$ which is taken by $k$ to a $1$-cube of $\mathcal{A}[\underline{n_2}]$. This implies that $\vd(k(x),k(y))=1$ by Proposition~\ref{newmetric} again. Consider a vertex $x$ of $\mathcal{A}[\underline{n_1}]$ and a sequence $(0=x_0,\dots,x_p=x)$ of vertices of $\mathcal{A}[\underline{n_1}]$ such that $\vd(x_i,x_{i+1})=1$ for all $i\in \{0,\dots,p-1\}$, which implies that $\vd(0,x)=p$ by Proposition~\ref{newmetric}. Then $\vd(k(x_i),k(x_{i+1}))=1$ for all $i\in \{0,\dots,p-1\}$, which implies $\vd(k(0),k(x))=p$ as well. From $\vd(0,n)=n$ by Proposition~\ref{newmetric}, we then deduce that $\vd(k(0),k(n))=n$. The only possibility is that $0=k(0)$ and that $n=k(n)$ by Proposition~\ref{newmetric}. We deduce 
\[
\forall x\in \mathcal{A}[\underline{n_1}]_0, \vd(0,x)=\vd(0,k(x)).
\]
Note that, in general, there is only the inequality $\vd(x,y)\geq \vd(k(x),k(y))$ for all vertices $x$ and $y$ of $\mathcal{A}[\underline{n_1}]$. For example, $x$ and $y$ can be incomparable vertices with $k(x)=k(y)$, which implies $\vd(x,y)=\infty$ and $\vd(k(x),k(y))=0$. The point is that $0$ is less than all vertices of $\mathcal{A}[\underline{n_1}]$ for the poset structure on vertices. This implies that for all  $p\in \vtx(\underline{n_2})\subset\vtx(\underline{n_1})$, $\vd(0,p)=\vd(0,k(p))=p$. We deduce \[\forall p\in \vtx(\underline{n_2}) \subset \vtx(\underline{n_1}), k(p)=p\] because there is only one vertex $x$ in the $\underline{n_1}$-cube and in the $\underline{n_2}$-cube such that $\vd(0,x)=p$. Thus the map \[k:\mathcal{A}[\underline{n_1}] \longrightarrow \mathcal{A}[\underline{n_2}]\] is of the form $k_1*\dots *k_{r}$ where \[k_j:\mathcal{A}[\underline{m_j}] \longrightarrow \mathcal{A}[|\underline{m_j}|]\] is a map of $\mathcal{A}$-sets for all $1\leq j \leq r$. We are then reduced to studying the case of a map of $\mathcal{A}$-sets \[f:\mathcal{A}[\underline{n}]\longrightarrow \mathcal{A}[n]\] with $\underline{n}=(n_1,\dots,n_p)$, $p\geq 2$ and $n=n_1+\dots+n_p$. By precomposing $f$ by the inclusion maps $\mathcal{A}[n_i]\subset \mathcal{A}[\underline{n}]$, we obtain a map of $\mathcal{A}$-sets 
\[
\begin{tikzcd}[row sep=3em, column sep=3em]
		f_i:\mathcal{A}[n_i]\arrow[r,"\subset"]& \mathcal{A}[\underline{n}] \arrow[r,"f"] &\mathcal{A}[n]
\end{tikzcd}
\]
for all $i\in\{1,\dots,p\}$. By Proposition~\ref{decomposition_distance_cor}, there exists a unique factorization 
\[
\begin{tikzcd}[row sep=3em, column sep=3em]
		\mathcal{A}[n_i]\arrow[d,equal] \arrow[rr,"f_i",bend left=30] \arrow[r,"\subset"] & \mathcal{A}[\underline{n}] \arrow[r,"f"] & \mathcal{A}[n] \arrow[d,equal]\\
\mathcal{A}[n_i] \arrow[r,"g_i"] & \mathcal{A}[n_i] \arrow[r,"\mathcal{L}_\mathcal{A}(h_i)"] & \mathcal{A}[n]
\end{tikzcd}
\]
for all $i\in\{1,\dots,p\}$. We obtain a factorization of $f:\mathcal{A}[\underline{n}]\to \mathcal{A}[n]$ as a composite 
\[
\begin{tikzcd}[row sep=3em, column sep=7em]
		\mathcal{A}[\underline{n}] \arrow[r,"\delta_{(g_1,\dots,g_p)}"] & \mathcal{A}[\underline{n}] \arrow[r,"\mathcal{L}_\mathcal{A}(h_1)*\dots *\mathcal{L}_\mathcal{A}(h_p)"] & \mathcal{A}[n]
\end{tikzcd}
\]
The right-hand map 
\[
\mathcal{L}_\mathcal{A}(h_1)*\dots *\mathcal{L}_\mathcal{A}(h_p):\mathcal{A}[\underline{n}] \to \mathcal{A}[n]
\]
is the image by the functor $\mathcal{L}_\mathcal{A}:\square^{op}\set\to \mathcal{A}^{op}\set$ of a map of precubical sets from $\square[\underline{n}]$ to $\square[n]$. It suffices to prove that it is a composite of maps of precubical sets of the form $\delta_{i,A,B}$ to complete the proof. There is nothing to prove when $\ell(\underline{n})= 1$. We make an induction on $\ell(\underline{n})\geq 2$. Assume first that $\underline{n}=(n_1,n_2)$. Consider a map of precubical sets $f:\square[\underline{n}] \to \square[n]$. By precomposing $f$ by the inclusion maps $\square[n_i]\subset \square[\underline{n}]$, we obtain a map of precubical sets \[f_i:\square[n_i]\subset \square[\underline{n}] \stackrel{f}\longrightarrow \square[n]\] for $i\in\{1,2\}$. Each map $f_i$ corresponds by the Yoneda lemma to an element $c_i$ of $\square[n]_{n_i}=\square([n_i],[n])$ such that the final state of $c_1$ is the initial state of $c_{2}$ and such that the initial state of $c_1$ is $0_n$ and such that the final state of $c_2$ is $1_n$. Thus there exists a partition $A_1\sqcup A_2=\{1,\dots,n\}$ such that $\de^0_{A_1}(c_1)=0_n$, $\de^1_{A_1}(c_1)=\de^0_{A_2}(c_2)$ and $\de^1_{A_2}(c_2)=1_n$. This implies that $c_1=\de^0_{A_2}(\id_{[n]})$ and $c_2=\de^1_{A_1}(\id_{[n]})$. We have proved that \[f=\phi_{A_1,A_2}\] which is the induction hypothesis for $p=2$. Consider now for some $p\geq 2$ a map of precubical sets 
\[
f:\square[\underline{n}] = (\square[n_1]*\dots *\square[n_p])*\square[n_{p+1}] \longrightarrow \square[n]
\]
with $n=n_1+\dots+n_{p+1}$ and $\underline{n}=(n_1,\dots,n_{p+1})$. The $\vd$ distance from the initial state of $\square[n_1]$ to the final state of $\square[n_p]$ in $\square[\underline{n}]$ being $n_1+\dots+n_p$, all cubes of $f(\square[n_1]*\dots *\square[n_p])$ are included in a subcube of $\square[n]$ of dimension $n_1+\dots+n_p$. This implies that the composite map \[\square[n_1]*\dots *\square[n_p]\subset \square[\underline{n}] \stackrel{f}\longrightarrow \square[n]\] factors as a composite of maps of precubical sets \[\square[n_1]*\dots *\square[n_p]\longrightarrow \square[n_1+\dots+n_p] \longrightarrow \square[n].\] We deduce that the map of precubical sets $f$ factors as a composite 
\[
f:(\square[n_1]*\dots *\square[n_p])*\square[n_{p+1}] \longrightarrow \square[n_1+\dots+n_p] *\square[n_{p+1}] \longrightarrow  \square[n].
\]
The induction hypothesis for $p\geq 2$ implies that the map of precubical sets $\square[n_1]*\dots *\square[n_p]\to \square[n_1+\dots+n_p]$ is a composite of maps of the form $\delta_{i,A,B}$. This implies that the left-hand map in the above factorization of $f$ is a composite of maps of the form $\delta_{i,A,B}$ as well. The induction hypothesis for $2$ implies that the right-hand map in the above factorization of $f$ is also a composite of maps of the form $\delta_{i,A,B}$. Hence the proof is complete.   
\epf

\begin{cor} \label{link}
	For any precubical $K$, Ziemia\'{n}ski's definition of $\Ch(K)$ given in \cite[Definition~1.1]{MR3722069} or \cite[Section~7]{MR4070250} and the above definition coincide.
\end{cor}

\bpf
Since $\square([n],[n])$ is a singleton for all $n\geq 0$, all maps of the form $\delta_{\underline{f}}$ are identities. 
\epf

\begin{nota}
	Since $\mathcal{L}_\mathcal{A}(\square[\underline{n}])=\mathcal{A}[\underline{n}]$ for all $\underline{n}\in \seq(n)$, $\mathcal{L}_\mathcal{A}$ being colimit-preserving, the functor $\mathcal{L}_\mathcal{A}:\square^{op}\set \longrightarrow \mathcal{A}^{op}\set$ induces a functor \[\mathcal{L}_\mathcal{A}^K:\Ch(K)\longrightarrow\Ch(\mathcal{L}_\mathcal{A}(K))\] for all precubical sets $K$.
\end{nota}

\bp \label{decomposition_distance_generalization}
Let $m\geq 0$. Let $K$ be a precubical set. Every map of $\mathcal{A}$-sets \[f:\mathcal{A}[m]\longrightarrow  \mathcal{L}_\mathcal{A}(K)\] factors uniquely as a composite of maps of $\mathcal{A}$-sets \[\mathcal{A}[m]\stackrel{g}\longrightarrow \mathcal{A}[m]\stackrel{\mathcal{L}_\mathcal{A}(h)}\longrightarrow \mathcal{L}_\mathcal{A}(K)\] where $h:\square[m]\to K$ is a map of precubical sets.
\ep

\bpf
The functor $L\mapsto L_m$ from $\mathcal{A}$-sets to sets is colimit-preserving. Thus, there are the bijections
\[
(\mathcal{L}_\mathcal{A}(K))_m \iso \liminj_{\square[p]\to K} (\mathcal{L}_\mathcal{A}(\square[p]))_m  \iso 
\liminj_{\square[p]\to K} (\mathcal{A}([p]))_m \iso 
\liminj_{\square[p]\to K} \mathcal{A}([m],[p]).
\]
The map of $\mathcal{A}$-sets $f:\mathcal{A}[m]\to  \mathcal{L}_\mathcal{A}(K)$ gives rise to an element of the set $(\mathcal{L}_\mathcal{A}(K))_m$ by Yoneda. Therefore there exists a map $g:\mathcal{A}[m]\to \mathcal{A}[p]$ and a map of precubical sets $h:\square[p]\to K$ such that $f=\mathcal{L}_\mathcal{A}(h)g$. From Proposition~\ref{decomposition_distance_cor} applied to $g$, we obtain $g=\mathcal{L}_\mathcal{A}(h')g'$ with $g'\in \mathcal{A}([m],[m])$ and we deduce $f=\mathcal{L}_\mathcal{A}(hh')g'$, which is the desired factorization. Consider two factorizations $f=\mathcal{L}(h_i)g_i$ with $h_i:\square[m]\to K$ and $g_i:\mathcal{A}[m]\to \mathcal{A}[m]$ with $i=1,2$. They correspond to two representatives of the map $f$ in the colimit of sets \[\liminj_{\square[p]\to K} \mathcal{A}([m],[p]).\] It means that there is a commutative diagram of $\mathcal{A}$-sets of the form
\[
\begin{tikzcd}[row sep=4em, column sep=4em]
	\mathcal{A}[m] \arrow[r,"g_1"] \arrow[d,equal] & \mathcal{A}[m] \arrow[d,"\mathcal{L}_\mathcal{A}(k_1)"] \arrow[r,"\mathcal{L}_\mathcal{A}(h_1)"]  & \mathcal{L}_\mathcal{A}(K) \arrow[d,equal]\\
\mathcal{A}[m]\arrow[d,equal] \arrow[r,"g_3"] & \mathcal{A}[p] \arrow[r,"\mathcal{L}_\mathcal{A}(h_3)"] & \mathcal{L}_\mathcal{A}(K)\arrow[d,equal]\\
\mathcal{A}[m] \arrow[r,"g_2"] & \mathcal{A}[m] \arrow[u,"\mathcal{L}_\mathcal{A}(k_2)"']\arrow[r,"\mathcal{L}_\mathcal{A}(h_2)"] & \mathcal{L}_\mathcal{A}(K)
\end{tikzcd}
\]
for some $g_3\in \mathcal{A}([m],[p])$ indexed by the map of precubical sets $h_3:\square[p]\to K$ in the diagram of sets. We obtain the equality $g_3 = \mathcal{L}_\mathcal{A}(k_1)g_1 = \mathcal{L}_\mathcal{A}(k_2)g_2$. From Proposition~\ref{decomposition_distance_cor}, we deduce $g_1=g_2$ and $k_1=k_2$. From the last equality, we deduce $\mathcal{L}_\mathcal{A}(h_1)=\mathcal{L}_\mathcal{A}(h_3)\mathcal{L}_\mathcal{A}(k_1)=\mathcal{L}_\mathcal{A}(h_3)\mathcal{L}_\mathcal{A}(k_2)=\mathcal{L}_\mathcal{A}(h_2)$. Since the functor $\mathcal{L}_\mathcal{A}:\square^{op}\set\to \mathcal{A}^{op}\set$ is faithful, we obtain $h_1=h_2$. Hence the proof is complete.
\epf

\bp \label{decomposition_distance_generalization2}
Let $\underline{n}\in \seq(n)$ with $n\geq 0$. Let $K$ be a precubical set. Every map of $\mathcal{A}$-sets $f:\mathcal{A}[\underline{n}]\longrightarrow  \mathcal{L}_\mathcal{A}(K)$ factors uniquely as a composite of maps of $\mathcal{A}$-sets 
\[
\begin{tikzcd}[row sep=small, column sep=3em]
	\mathcal{A}[\underline{n}]\arrow[r,"g"] & \mathcal{A}[\underline{n}]\arrow[r,"\mathcal{L}_\mathcal{A}(h)"] &\mathcal{L}_\mathcal{A}(K)
\end{tikzcd}
\]
where $h:\square[\underline{n}]\to K$ is a map of precubical sets. 
\ep

\bpf
Let $\underline{n}=(n_1,\dots,n_p) \in \seq(n)$. Giving $f$ is equivalent to giving $p$ maps of $\mathcal{A}$-sets $f_i:\mathcal{A}[n_i]\to \mathcal{L}_\mathcal{A}(K)$ for $1\leq i \leq p$ satisfying $f_i(1_{n_i}) = f_{i+1}(0_{n_{i+1}})$ for $1\leq i\leq p-1$. Write \[f=f_1*\dots * f_p.\] In the same way, given $h$ is equivalent to giving $p$ maps of precubical sets $h_i:\square[n_i]\to K$ for $1\leq i \leq p$ satisfying $h_i(1_{n_i}) = h_{i+1}(0_{n_{i+1}})$ for $1\leq i\leq p-1$. Write 
\[
h=h_1*\dots *h_p \hbox{ and }\mathcal{L}_\mathcal{A}(h)=\mathcal{L}_\mathcal{A}(h_1)*\dots *\mathcal{L}_\mathcal{A}(h_p).
\]
Then $f=\mathcal{L}_\mathcal{A}(h)g$ if and only if $f_i=\mathcal{L}_\mathcal{A}(h_i)g_i$ for all $1\leq i\leq p$ where $g_i:[n]\to [n]$ is a map of $\mathcal{A}$. The proof is complete thanks to Proposition~\ref{decomposition_distance_generalization}. 
\epf

\begin{figure}
\[
\begin{tikzcd}[row sep=3em, column sep=2em]
\mathcal{A}[\underline{m}]	\arrow[rd,equal] \arrow[dd,"\underline{c}"'] \arrow[rr,"g_0"] && \mathcal{L}_\mathcal{A}(\square[\underline{m}])
\arrow[rd,dashed,"\mathcal{L}_\mathcal{A}(h_1)"] \arrow[dd,"\mathcal{L}_\mathcal{A}(h_0)",pos=0.8]  &\\
& \mathcal{A}[\underline{m}] \arrow[rr,"g",pos=0.3,crossing over]   && \mathcal{L}_\mathcal{A}(\square[\underline{n}]) \arrow[dd,"\mathcal{L}_\mathcal{A}(h)"] \\
\mathcal{L}_\mathcal{A}(K) \arrow[rd,equal] \arrow[rr,equal]  &&
\mathcal{L}_\mathcal{A}(K) \arrow[rd,equal]& \\
& \arrow[from=uu,"\underline{c}"',pos=0.3,crossing over]\mathcal{L}_\mathcal{A}(K) \arrow[rr,equal] &&
\mathcal{L}_\mathcal{A}(K)
\end{tikzcd}
\]
\caption{Construction of $h_1$}
\label{h1}
\end{figure}

\bth \label{hom-eq-cube-chain}
Let $K$ be a precubical set. The functor \[\mathcal{L}_\mathcal{A}^K:\Ch(K)\to\Ch(\mathcal{L}_\mathcal{A}(K))\] induces a homotopy equivalence \[|\Ch_{\alpha,\beta}(K)|\simeq|\Ch_{\alpha,\beta}(\mathcal{L}_\mathcal{A}(K))|\] for all $(\alpha,\beta)\in K_0\p K_0$ between the classifying spaces where $|\C|$ means the classifying space of $\C$, i.e. the geometric realization of the simplicial nerve of $\C$.
\eth

\bpf
By Quillen's Theorem~A \cite[Theorem~19.6.14]{ref_model2}, it suffices to prove that the comma category $(\underline{c}\ddownarrow \mathcal{L}_\mathcal{A}^K)$ is contractible for all cube chains $\underline{c}$ of $\mathcal{L}_\mathcal{A}(K)$. By \cite[Proposition~14.3.14]{ref_model2}, it suffices to prove that the comma category $(\underline{c}\ddownarrow \mathcal{L}_\mathcal{A}^K)$ has an initial object for all cube chains $\underline{c}$ of $\mathcal{L}_\mathcal{A}(K)$. An object of the comma category $(\underline{c}\ddownarrow \mathcal{L}_\mathcal{A}^K)$ is a commutative square of $\mathcal{A}$-sets of the form 
\[\begin{tikzcd}[row sep=3em, column sep=3em]
	\mathcal{A}[\underline{m}]\arrow[d,"\underline{c}"'] \arrow[r,"g"] & \mathcal{L}_\mathcal{A}(\square[\underline{n}]) \arrow[d,"\mathcal{L}_\mathcal{A}(h)"]\\
\mathcal{L}_\mathcal{A}(K) \arrow[r,equal] &   \mathcal{L}_\mathcal{A}(K)
\end{tikzcd}\]
such that $|\underline{m}| = |\underline{n}|$ and $\vtx(\underline{n})\subset \vtx(\underline{m})$. Proposition~\ref{decomposition_distance_generalization2} provides a commutative square of $\mathcal{A}^{op}\set$
\[
\begin{tikzcd}[row sep=3em, column sep=3em]
	\mathcal{A}[\underline{m}]\arrow[d,"\underline{c}"'] \arrow[r,"g_0"] & \mathcal{L}_\mathcal{A}(\square[\underline{m}]) \arrow[d,"\mathcal{L}_\mathcal{A}(h_0)"]\\
\mathcal{L}_\mathcal{A}(K) \arrow[r,equal] &   \mathcal{L}_\mathcal{A}(K).
\end{tikzcd}
\]
which is clearly an object of the comma category $(\underline{c}\ddownarrow \mathcal{L}_\mathcal{A}^K)$. Consider the diagram of solid arrows of $\mathcal{A}$-sets depicted in Figure~\ref{h1} where the front face is another object of the comma category $(\underline{c}\ddownarrow \mathcal{L}_\mathcal{A}^K)$. This implies that $|\underline{m}| = |\underline{n}|$ and $\vtx(\underline{n}) \subset \vtx(\underline{m})$. It can be rearranged as follows:
\[
\begin{tikzcd}[row sep=4em, column sep=5em]
	\mathcal{A}[\underline{m}] \arrow[d,equal] \arrow[r,bend left=10pt,"g_0"] \arrow[r,dashed,bend right=10pt,"g_1"']  & \mathcal{L}_\mathcal{A}(\square[\underline{m}]) \arrow[d,dashed,"\mathcal{L}_\mathcal{A}(h_1)"] \arrow[r,"\mathcal{L}_\mathcal{A}(h_0)"] & \mathcal{L}_\mathcal{A}(K) \arrow[d,equal]\\ 
\mathcal{A}[\underline{m}] \arrow[r,"g"] & \mathcal{L}_\mathcal{A}(\square[\underline{n}]) \arrow[r,"\mathcal{L}_\mathcal{A}(h)"] & \mathcal{L}_\mathcal{A}(K)
\end{tikzcd}
\]
with $\underline{c}=\mathcal{L}_\mathcal{A}(h_0)g_0=\mathcal{L}_\mathcal{A}(h)g$. Using Proposition~\ref{decomposition_distance_generalization2}, write $g=\mathcal{L}_\mathcal{A}(h_1)g_1$. We obtain the equalities \[\mathcal{L}_\mathcal{A}(h_0)g_0 =\underline{c}=\mathcal{L}_\mathcal{A}(h)g=\mathcal{L}_\mathcal{A}(h)\mathcal{L}_\mathcal{A}(h_1)g_1=\mathcal{L}_\mathcal{A}(hh_1)g_1.\] By Proposition~\ref{decomposition_distance_generalization2} and since $\mathcal{L}_\mathcal{A}$ is faithful by Proposition~\ref{LA-faithful}, we obtain \[g_1=g_0 \hbox{ and }hh_1=h_0.\] We have proved that there exists a diagram of solid arrows of $\mathcal{A}$-sets
\[
\begin{tikzcd}[row sep=5em, column sep=5em]
	\mathcal{A}[\underline{m}] \arrow[d,equal] \arrow[r,"g_0"]  & \mathcal{L}_\mathcal{A}(\square[\underline{m}]) \arrow[d,"\mathcal{L}_\mathcal{A}(h_1)"] \arrow[r,"\mathcal{L}_\mathcal{A}(h_0)"] & \mathcal{L}_\mathcal{A}(K) \arrow[d,equal]\\ 
\mathcal{A}[\underline{m}] \arrow[r,"g"] & \mathcal{L}_\mathcal{A}(\square[\underline{n}]) \arrow[r,"\mathcal{L}_\mathcal{A}(h)"] & \mathcal{L}_\mathcal{A}(K)
\end{tikzcd}
\]
From $g=\mathcal{L}_\mathcal{A}(h_1)g_0$ and Proposition~\ref{decomposition_distance_generalization2}, we deduce that there exists a unique map $\mathcal{L}_\mathcal{A}(h_1)$ making the above diagram commutative, and therefore a unique $h_1$, the functor $\mathcal{L}_\mathcal{A}$ being faithful by Proposition~\ref{LA-faithful}. Since $|\underline{m}| = |\underline{n}|$ and $\vtx(\underline{n}) \subset \vtx(\underline{m})$ by hypothesis, the commutative square of precubical sets 
\[
\begin{tikzcd}[row sep=3em, column sep=3em]
	\square[\underline{m}] \arrow[d,"h_0"'] \arrow[r,"h_1"] & \square[\underline{n}] \arrow[d,"h"]\\
K \arrow[r,equal] &   K
\end{tikzcd}
\]
is a map of $\Ch_{\alpha,\beta}(K)$. Thus the diagram of solid arrows depicted in Figure~\ref{h1} yields together with $h_1$ a well-defined map of the comma category $(\underline{c}\ddownarrow \mathcal{L}_\mathcal{A}^K)$. This implies that the commutative square 
\[
\begin{tikzcd}[row sep=3em, column sep=3em]
	\mathcal{A}[\underline{m}]\arrow[d,"\underline{c}"'] \arrow[r,"g_0"] & \mathcal{L}_\mathcal{A}(\square[\underline{m}]) \arrow[d,"\mathcal{L}_\mathcal{A}(h_0)"]\\
\mathcal{L}_\mathcal{A}(K) \arrow[r,equal] &   \mathcal{L}_\mathcal{A}(K)
\end{tikzcd}
\]
is an initial object of the comma category $(\underline{c}\ddownarrow \mathcal{L}_\mathcal{A}^K)$ and the proof is complete.
\epf

We recall \cite[Definition~4.11]{DirectedDegeneracy} which is necessary to understand Corollary~\ref{Zgen}. There are the homeomorphisms (natural with respect to $[n]\in \square$)
\[
|\square[n]|_{geom} \iso [0,1]^n \iso |\mathcal{A}[n]|_{geom}
\]
by definition of the geometric realization. A \textit{tame $d$-path} of $|\mathcal{A}[n]|_{geom} = [0,1]^n$ is a nonconstant continuous map $\gamma:[0,\ell]\to [0,1]^n$ with $\ell>0$ such that $\gamma(0),\gamma(\ell)\in \{0,1\}^n$ and such that $\gamma$ is nondecreasing with respect to each axis of coordinates. Let $c\in K_n$ with $n\geq 1$ be an $n$-cube of an $\mathcal{A}$-set $K$. A \textit{tame $d$-path} of $c$ is a composite continuous map denoted by $[c;\gamma]:[0,\ell] \to |K|_{geom}$ with $\ell>0$ such that $\gamma:[0,\ell]\to [0,1]^n$ is a $d$-path with $[c;\gamma]=|c|_{geom}\gamma$. Let $K$ be a general $\mathcal{A}$-set. A \textit{tame $d$-path} of $K$ is a continuous path $[0,\ell] \to |K|_{geom}$ which is a Moore composition $[c_1;\gamma_1] * \dots *[c_n;\gamma_n]$ of $d$-paths of the cubes $c_1,\dots,c_n$ of $K$. $\gamma(0)\in K_0$ is called the \textit{initial state} of $\gamma$ and $\gamma(\ell)\in K_0$ is called the \textit{final state} of $\gamma$.

\begin{cor} \label{Zgen} 
	Let $K$ be a precubical set. The space of tame natural $d$-paths of $\mathcal{L}_\mathcal{A}(K)$ is homotopy equivalent to $|\Ch(\mathcal{L}_\mathcal{A}(K))|$.
\end{cor}

\bpf
There are the homeomorphisms (natural with respect to $[n]\in \square$)
\[
|\square[n]|_{geom} \iso [0,1]^n \iso |\mathcal{A}[n]|_{geom} \iso |\mathcal{L}_\mathcal{A}(\square[n])|_{geom},
\]
by Proposition~\ref{free_square}. Since all involved functors are colimit-preserving, we obtain for all precubical sets $K$ the natural homeomorphism 
\[
|K|_{geom}\iso |\mathcal{L}_\mathcal{A}(K)|_{geom}.
\]
Hence the space of tame natural $d$-paths of $K$ is equal to the one of $\mathcal{L}_\mathcal{A}(K)$. By \cite[Theorem~7.5]{MR4070250}, we deduce that the space of tame natural $d$-paths of $\mathcal{L}_\mathcal{A}(K)$ is homotopy equivalent to $|\Ch(K)|$, and therefore homotopy equivalent to $|\Ch(\mathcal{L}_\mathcal{A}(K))|$ by Theorem~\ref{hom-eq-cube-chain}. 
\epf

\section{Application}
\label{app}

We refer to \cite{ccsprecub,symcub} for further details. \cite[Theorem~4.1.8]{symcub} states that, for two labelled cubes $\square[a_1,\dots,a_m]$ and $\square[a_{m+1},\dots,a_{m+n}]$ with $m\geq 0$ and $n\geq 0$, there is an isomorphism of $\sigma$-labelled symmetric transverse sets
\begin{multline*}
	\mathcal{L}_{\widehat{\square}_S}\big(\COSK^\Sigma( \square[a_1,\dots,a_m]_{\leq 1}
	\p_\Sigma \square[a_{m+1},\dots,a_{m+n}]_{\leq 1}) \big)\\\iso
	\cosk_1^{\widehat{\square}_S,\overline{\Sigma}}\big(
	\widehat{\square}_S[a_1,\dots,a_m]_{\leq 1} \overline{\p}_\Sigma
	\widehat{\square}_S[a_{m+1},\dots,a_{m+n}]_{\leq 1}\big)
\end{multline*}
where 
\begin{itemize}[leftmargin=*]
	\item $\COSK^\Sigma$ is the $\sigma$-labelled directed coskeleton of \cite[Section~3.3]{ccsprecub} which is a tweak of the $\sigma$-labelled coskeleton functor of the category of precubical sets, the latter being badly behaved by \cite[Proposition~3.15]{ccsprecub}: it contains too many cubes and some of them have to be identified;
	\item $\cosk_1^{\widehat{\square}_S,\overline{\Sigma}}$ is the $\sigma$-labelled coskeleton functor of the category of symmetric transverse sets which is well behaved by \cite[Theorem~3.1.15]{symcub}.
\end{itemize}
By Theorem~\ref{hom-eq-cube-chain}, the underlying precubical set of \[\COSK^\Sigma( \square[a_1,\dots,a_m]_{\leq 1}
\p_\Sigma \square[a_{m+1},\dots,a_{m+n}]_{\leq 1})\] and the underlying symmetric transverse set of \[\cosk_1^{\widehat{\square}_S,\overline{\Sigma}}\big(
\widehat{\square}_S[a_1,\dots,a_m]_{\leq 1} \overline{\p}_\Sigma
\widehat{\square}_S[a_{m+1},\dots,a_{m+n}]_{\leq 1}\big)\] have homotopy equivalent categories of cube chains. 

This implies that the parallel composition with synchronization for process algebra, as it is formalized in \cite{symcub} using the labelled coskeleton functor of the category of symmetric transverse sets, has a category of cube chains which gives the correct homotopy type of tame natural $d$-paths.

\end{document}